\documentclass[12pt]{article}
\usepackage[english]{babel}
\usepackage[cp1251]{inputenc}
\usepackage{amsthm}

\begin{document}
\newtheorem{defn}{Definition}[section]
\newtheorem{lm}{Lemma}[section]
\newtheorem{thm}{Theorem}[section]
\newtheorem{pr}{Proposition}[section]
\begin{center}
{\bf Conjugacy of Cartan subalgebras of complex finite
dimensional Leibniz algebras}\\[3mm]
{\sl Omirov B.A.}\footnote{The work is supported by grant of the
Academy of Sciences of the Republic of Uzbekistan No. 39-04 and
grant of INTAS Ref. Nr. 04-83-3035}
\end{center}
\begin{abstract}

In the present work the properties of Cartan subalgebras and their
connection with regular elements in finite dimensional Lie
algebras are extended to the case of Leibniz algebras. It is shown
that Cartan subalgebras and regular elements of a Leibniz algebra
correspond to Cartan subalgebras and regular elements of a Lie
algebra by a natural homomorphism. Conjugacy of Cartan subalgebras
of Leibniz algebras is shown.

\end{abstract}

\section{Introduction}

This work is devoted to the study of Leibniz algebras, which were
introduced by French mathematician J.-L. Loday in \cite{Lod} and
considered further in works \cite{Bal} - \cite{Gao}.

It is well known that Leibniz algebras are "non commutative"
generalization of Lie algebras. Investigations of nilpotent
Leibniz algebras in \cite{AyupOmir1} - \cite{AlAyupOmir2} show
that many nilpotent properties of Lie algebras can be extended to
the case of nilpotent Leibniz algebras.

A characteristic of non-Lie Leibniz algebras is a nontriviality of
the ideal generated by squares of the elements of the algebra
(moreover, it is abelian).

Cartan subalgebras, their relations with regular elements and
decomposition to weight spaces by Cartan subalgebras play the
basic role in the structure theory of Lie algebras. Although some
properties of regular elements and Cartan subalgebras were studied
in \cite{AlAyupOmir2}, the relation between regular elements and
Cartan subalgebras of the Leibniz algebra, on the one side, and
regular elements and Cartan subalgebras of its factor algebra by
the ideal generated by squares of elements (which evidently is a
Lie algebra), on the other side, has been not clarified.
Therefore, one of the aim of this work is to investigate the
mentioned relation. Namely, it is proved that the images of
regular elements and Cartan subalgebras by a natural homomorphism
are regular elements and Cartan subalgebras respectively.

The classical result of the structure theory of finite dimensional
Lie algebras on conjugacy of Cartan subalgebras is proven in the
case of Leibniz algebras.

\section{Preliminaries}

\begin{defn}
An algebra $L$ over a field $F$ is called a Leibniz algebra if the
Leibniz identity
$$
[x,[y,z]]=[[x,y],z]-[[x,z],y]
$$
holds for any $x,y,z\in L,$ where $[\ ,\ ]$ is the multiplication
in $L.$ \end{defn}

Note that in the case of fulfilling the identity $[x,x]=0$ the
Leibniz identity can be easily reduced to the Jacobi identity.
Thus, Leibniz algebras are "non commutative" analogues of Lie
algebras.

For an arbitrary algebra $L,$ we define the following sequence:
$$
L^1=L,\qquad L^{n+1}=[L^n, L^1].
$$
\begin{defn}
An algebra $L$ is called nilpotent if there exists $s\in N$ such
that $L^s=0.$
\end{defn}

For an arbitrary element $x\in L,$ we consider the operator of
right multiplication $R_x:L\to L$ where $R_x(z)=[z,x].$ The set
$R(L)=\{ R_x: x\in L\}$ is a Lie algebra with respect to operation
of commutation, and the following identity holds:
$$
R_xR_y-R_yR_x=R_{[y,x]}.
$$

It is easy to see from this identity that the solvability of the
Lie algebra $R(L)$ is equivalent to the solvability of the Leibniz
algebra $L.$ Moreover, if $L$ is nilpotent then the Lie algebra
$R(L)$ is also nilpotent.

The following lemma gives a decomposition of a vector space into
the direct sum of two invariant subspaces with respect to a linear
transformation.

\begin{lm}[Fitting's Lemma] Let $V$ be a vector space and
$A:V\to V$ be a linear transformation. Then $V= V_{0A}\oplus
V_{1A}$, where $A(V_{0A})\subseteq V_{0A},$ $A(V_{1A})\subseteq
V_{1A}$ and $V_{0A}=\{ v\in V|\ A^i(v)=0$ for some $i\}$ and
$V_{1A}=\bigcap\limits_{i=1}^\infty A^i(V).$ Moreover, $A_{|
V_{0A}}$ is a nilpotent transformation and $A_{V_{1A}}$ is an
automorphism. \end{lm}

\begin{proof}
See \cite{Jac} (chapter II, \S 4).
\end{proof}

\begin{defn} The spaces $V_{0A}$ and $V_{1A}$ are called
the Fitting's null-component and the Fitting's one-component
(respectively) of the space $V$ with respect to the transformation
$A.$ \end{defn}

\begin{defn} An element $h$ of the Leibniz algebra $L$ is said to
be regular if the dimension of the Fitting's null-component of the
space $L$ with respect to $R_h$ is minimal. In addition, its
dimension is called a rank of the algebra $L.$
\end{defn}

It is easy to see that the dimension of the Fitting's
null-component of a linear transformation $A$ equals to the order
of zero root of characteristic polynomial of this transformation.
Hence an element $h$ is regular if and only if the order of zero
characteristic root is minimal for $R_h.$

Note that in the case of Lie algebras the linear transformation
$R_h$ is degenerated (since $[h,h]=0$ for any $h$) and therefore
the rank of the Lie algebra is greater than zero.

The following lemma shows that for Leibniz algebras the rank is
also greater than zero.

\begin{lm} Let $L$ be a finite dimensional Leibniz algebra. Then
the operator $R_x$ is degenerated for any $x\in L.$ \end{lm}
\begin{proof}

See \cite{AlAyupOmir2}.
\end{proof}

We also have the following generalization of Fitting's Lemma for
Lie algebras of nilpotent transformations of a vector space.

\begin{thm} Let $G$ be a nilpotent Lie algebra of linear
transformations of a vector space $V$ and $V_0=\bigcap\limits_{A
\in G} V_{0A},$ $V_1=\bigcap\limits_{ i=1}^\infty G^i(V).$ Then
the subspaces $V_0$ and $V_1$ are invariant with respect to $G$
(i.e. they are invariant with respect to every transformation $B$
of $G$) and $V=V_0\oplus V_1.$ Moreover, $V_1=\sum\limits_{A\in G}
V_{1A}.$ \end{thm}
\begin{proof}
See \cite{Jac} (chapter II, \S 4).
\end{proof}

\emph{Remark 1.} From \cite{Jac} in case of vector space V over an
infinite field and fulfillment of conditions of theorem 2.1, we
have the existence of the element $B\in G$ such that $V_0= V_{0B}$
and $V_1=V_{1B}$.

\section{Conjugacy of the Cartan subalgebras of finite dimensional
Leibniz algebras.}

Let $\Im$ be a nilpotent subalgebra of a Leibniz algebra $L$ and
$L=L_0\oplus L_1$ be the Fitting's decomposition of the algebra
$L$ with respect to the nilpotent Lie algebra $R(\Im)=\{ R_x|\
x\in \Im \}$ of transformations of the vector space as in theorem
2.1.

The set $l(\Im)=\{ x\in L|\ [x,\Im]\subseteq \Im\}$ is said to be
a left normalizator of the subalgebra $\Im$ in the algebra $L.$

The set $r(\Im)=\{ x\in L|\ [\Im,x]\subseteq \Im\}$ is said to be
a right normalizator of the subalgebra $\Im$ in the algebra $L.$

\begin{defn} A subalgebra $\Im$ of a Leibniz algebra $L$ is called
a Cartan subalgebra if the following two conditions are satisfied:\\
\indent a) $\Im$ is nilpotent;\\
\indent b) $\Im$ coincides with the left normalizator of $\Im$ in
the algebra $L.$ \end{defn}

Note that the definition of a Cartan subalgebra of a Leibniz
algebra is agree with the definition of a Cartan subalgebra for
the Lie algebra.

In view of the antisymmetric identity in Lie algebras, the sets
$l(\Im)$ and $r(\Im)$ coincide.

The following example show that in general these sets do not
coincide for Leibniz algebras.

\emph{\bf Example 3.1.} Let $L$ be a Leibniz algebra determined by
the following multiplication:
$$
[x,z]=x,\quad [z,y]=y,\quad [y,z]=-y,\quad [z,z]=x,
$$
where $\{ x,y,z\}$ is the basis of the algebra $L$ and omitted
products are equal to zero.\\
\indent Then $\Im=\{ x-z\}$ is a Cartan subalgebra, but $r(\Im)=\{
x,z\}.$

For Cartan subalgebras of Leibniz algebras similar to the case of
Lie algebras, there is a characterization in terms of the
Fitting's null-component, namely, the following proposition is
true.
\begin{pr} A nilpotent subalgebra $\Im$ of a Leibniz algebra $L$
is a Cartan subalgebra if and only if $\Im$ coincides with $L_0$
in the Fitting's decomposition of the algebra $L$ with respect to
$R(\Im).$ \end{pr}
\begin{proof} See \cite{AlAyupOmir2}. \end{proof}

The following theorem establishes properties of the Fitting's
null-component of the regular element of a Leibniz algebra.
\begin{thm} Let $L$ be a Leibniz algebra over an infinite field $F$
and $a$ be a regular element of the algebra $L.$ Then the
Fitting's null-component $\Im$ of the algebra $L$ with respect to
$R_a$ is a Cartan subalgebra. \end{thm}

\begin{proof} See [10]. \end{proof}

Another useful remark on regular elements and Cartan subalgebras
of Leibniz algebras is the fact that if the Cartan subalgebra
contains a regular element $a,$ then $\Im$ is uniquely determined
by the element $a$ as the Fitting's null-component of the algebra
$L$ with respect to $R_a,$ i.e. $\Im=L_{0R_a}.$

For the Leibniz algebra $L,$ we consider the natural homomorphism
$\varphi$ into the factor algebra $L/I,$ where $I=ideal\langle
[x,x]|\ x\in L\rangle.$

\begin{pr} Let $L$ be a complex finite dimensional Leibniz
algebra.Then the image of a regular element of the algebra $L$ by
a homomorphism $\varphi$ is a regular element of the Lie algebra
$L/I.$ \end{pr}
\begin{proof} Let $a$ be a regular element of the algebra $L.$
We will prove that the element $\overline a=a+I$ will be a regular
element of the Lie algebra $L/I.$ Suppose the opposite, i.e.
$\overline a=a+I$ is not a regular element. Let $\overline b=b+I$
be any regular element of the Lie algebra $L/I$ and $a-b\not\in
I.$

Since $I$ is an ideal, then for any $x\in L$ we have
$R_x(I)\subseteq I.$ It means that the matrix of the
transformation $R_x$ has the following block form
$$ R_x=\left(\begin{array}{cc} X,& 0 \\ Z_X, & I_X \\ \end{array}
\right) $$ in the basis $\{ e_1,e_2,\ldots, e_m,i_1,i_2,\ldots,
i_n\}$ of $L,$ where $\{ i_1,i_2,\ldots, i_n\}$ is the basis of
$I.$ Here $X$ is the matrix of the transformation $R_x|_{ \{
e_1,e_2,\ldots, e_m\}}$ and $I_X$ is the matrix of the
transformation $R_x|_I.$

Let
$$
R_a=\left(\begin{array}{cc} A,& 0 \\ Z_a, & I_a \\ \end{array}
\right),\quad R_b=\left(\begin{array}{cc} B,& 0 \\ Z_b, & I_b \\
\end{array} \right)
$$
be the matrices of the transformations $R_a$ and $R_b$
respectively.

Let $k$ (respectively $k'$) be the order of the characteristic
zero root of the matrix $A$ (respectively $B$) and $s$ and $s'$ be
the orders of the characteristic zero root of the matrices $I_a$
and $I_b,$ respectively. Then we have $k'<k,$ $s<s'.$

Put $U=\bigg\{ y\in L\setminus I \left|\ R_y=\left(\begin{array}{cc} Y,& 0 \\
Z_y, & I_y \\ \end{array} \right)\right.$ and $Y$ has the order of
the characteristic zero root less than $k \bigg\},$ $V=\bigg\{
y\in L\setminus I \left|\ R_y=\left(\begin{array}{cc} Y,& 0 \\
Z_y, & I_y \\ \end{array} \right)\right.$ and $I_y$ has the order
of the characteristic zero root less than $s+1 \bigg\}.$

Since $b\in U$ and $a\in V$ these sets are non empty.

Let's show that the set $U$ is an open subset of the set
$L\setminus I$ in the Zariski topology.

Let $Y$ have the order of the characteristic zero root less than
$k.$ Then $Y^k$ has the rank greater than $n-k.$ It means that
there exists a non-zero minor of the order $n-k+1.$ In the other
words, there exists a non-zero polynomial of structural constants
of the algebra $L$, hence the set $U$ is open in the Zariski
topology in the subset of the set $L\setminus I.$

One can analogously prove that the set $V$ is open in $L\setminus
I.$ It is not difficult to check that the sets $U$ and $V$ are
dense in $L\setminus I.$ Therefore, there exists an element $y\in
U\cap V$ such that $Y$ has the order of characteristic zero root
less than $k$ and $I_y$ has the order of the characteristic zero
root less than $s+1.$ Thus, for this element $y$ the order of
characteristic zero root is not greater than $k+s-1,$ i.e. the
rank of the algebra $L$ is less than $k+s$ and we obtain the
contradiction to the assumption that $\overline a$ is not a
regular element of the Lie algebra $L\setminus I.$
\end{proof}

\emph{Remark 2.} For the Cartan subalgebra $\Im$ of the Leibniz
algebra $L,$ we consider the Lie algebra $R(\Im)$ of linear
transformations $L$ (which evidently is nilpotent) and the
decomposition of the algebra $L$ with respect to $R(\Im).$ Remark
1 implies existence of an element $R_b\in R(\Im)$ such that the
Fitting's null-component with respect to the nilpotent Lie algebra
of linear transformations $R(\Im)$ coincides with the Fitting's
null component with respect to the transformation $R_b,$ i.e.
$L_0= L_{0R_b}.$ Using Proposition 3.1 we obtain $\Im= L_{0R_b}.$

Let $\Im$ be a Cartan subalgebra of the Leibniz algebra $L$ and
$L=L_0\oplus L_\alpha\oplus L_\beta\oplus \ldots\oplus L_\gamma$
be a decomposition of the algebra $L$ into characteristic
subspaces with respect to operator $R_b$ possessing the property
$\Im=L_{0R_b}.$
\begin{lm} Let the element
$x=x_0+x_\alpha+x_\beta+\ldots+x_\gamma,$ where $x_\sigma\in
L_\sigma,$ $\sigma\in\{0,\alpha ,\beta, \ldots, \gamma\}$ satisfy the following conditions: \\
\indent a) there exists $k\in N$ such that $[\ldots[[
x,\underbrace{b],b],\ldots,b}_{\mbox{$k$
times}}]\in I $;\\ \indent b) $x\ne x_0.$ \\
\indent Then $\overline x=\overline x_0.$ \end{lm}
\begin{proof} Let $k$ is the minimal of numbers having
the property: $[\ldots[[ x,\underbrace{b],b],\ldots,b}_{\mbox{$k$
times}}]\in I.$ It is not difficult to see that the element
$[\ldots[[ x,\underbrace{b],b],\ldots,b}_{\mbox{$k-1$ times}}]\in
I $ possesses the properties of the element $x.$ That is why we
can suppose $[x,b]\in I.$

To prove the lemma, it is sufficient to show that
$x'=x_\alpha+x_\beta+\ldots+x_\gamma\in I.$ Without any loss of
generality we can also suppose that $[x',b]\in I.$

If $x_\sigma\in I$ for any $\sigma\in\{\alpha, \beta, \ldots,
\gamma\}$ than the assertion of the lemma is evident. Suppose
there exists $\sigma\in\{\alpha, \beta, \ldots, \gamma\},$ such
that $x _\sigma\not\in I.$ For convenience, we put
$\sigma=\alpha.$

Let $x_\alpha=\alpha_1e_1+\alpha_2e_2+\ldots \alpha_te_t$ be the
decomposition with respect to Jordan basis of $L_\alpha.$ Then
$[x_\alpha,b]=\alpha x_\alpha+\alpha_1'e_2+\alpha_2'e_3+\ldots+
\alpha_{t-1}'e_t$ for some $\alpha_i'$, where $i=\{ 1,2,\ldots,
t-1 \}.$

Consider multiplication:
$$
[x',b] =\alpha x'+\alpha_1'e_2+\alpha_2' e_3+ \ldots+
\alpha_{t-1}'e_t +\left( [x_\beta,b]-\alpha x_\beta\right)+\ldots+
\left( [x_\gamma,b]-\alpha x_\gamma\right).
$$

For the element $x''=[x',b]-\alpha x'$ (which has not the basis
element $e_1$ in its decomposition), we have $[x'',b]\in I$ and
$\overline{x''}=\alpha\overline{x'}$, and also if $x''\in L_0$
then $x''=0.$ Hence, $x'=x_\alpha$ and $[x',b]=\alpha x'$ which
implies that $x_\alpha\in I$. This contradicts the assumption that
$x_\alpha\not\in I.$

Thus, conditions $[x'',b]\in I;$ $x''\ne x_0$ and $\overline{x''}
=\alpha\overline{x'}$ hold for the element $x''.$

Let $x''=x_\alpha'+x_\beta'+\ldots+x_\gamma',$ where the basis
element $e_1$ is absent in the decomposition of $x_\alpha'.$ If
$x_\alpha'\in I$ than we can suppose that $x''$ equals to
$x_\beta'+\ldots+x_\gamma'.$

If $x_\alpha'\not\in I$ then we analogously obtain the element
$x'''=x_\alpha''+x_\beta''+\ldots+x_\gamma''$ satisfying the
following conditions: $[x''',b]\in I,$ $\overline{x'''}
=\alpha^2\overline{x''}$ and there are no basis elements $e_1,$
$e_2$ in the decomposition of $x_\alpha''.$

If we continue similar operations finite times and take into
account the property of the Jordan basis, we obtain existence of
the element $z=z_\alpha+z_\beta+\ldots+z_\gamma$ such that
$[z_\alpha ,b]=\alpha z_\alpha$ and conditions $[z,b]\in I,$ $z\ne
x_0$ and $\overline z=\alpha^q\overline{x'}$ are satisfied.

Consider the element $z'=[z,b]-\alpha z=([z_\beta,b]-\alpha
z_\beta)+\ldots+([z_\gamma,b]-\alpha z_\gamma).$ It is evident
that $[z',b]\in I$ and $\overline{z'}=\alpha^{q+1}\overline{x'}.$
In the case of $z'=x_0,$ we have $z_\beta=\ldots=z_\gamma=0,$ i.e.
$z'=0$ and $x'\in I.$

Suppose that $z'\ne x_0.$ In this case we have shown the existence
of the element $z'=z_\beta+\ldots+z_\gamma$ (such that there is
not a component lying in $L_\alpha$ in the decomposition of $z'$)
for which the following conditions hold: $[z',b]\in I;$ $z\ne x_0$
and $\overline{z'}=\alpha^{q+1}\overline{x'}.$ If we continue
analogously, we obtain the existence of an element $d$ such that
$[d,b]=\delta d$ and $d\in L_\delta$, where $\delta\in\{\alpha,
\beta,\ldots,\gamma\},$ satisfying the following conditions:
$[d,b]\in I,$ $\overline{d}=q\overline{x'}$ for a non-zero number
$q.$ Hence, $d\in I$ and therefore $\overline{x'}= \overline{0},$
i.e. $x'\in I.$
\end{proof}

The following theorem establishes a relation for Cartan
subalgebras as in Proposition 3.2.
\begin{thm} Let $L$ be a Leibniz algebra and $\Im$ is its Cartan
subalgebra. Then the image of the subalgebra $\Im$ by homomorphism
$\varphi$ is the Cartan subalgebra of the Lie algebra $L/I.$
\end{thm}
\begin{proof} Denote $\varphi(\Im)= \overline \Im.$ Nilpotence of $\overline \Im$ follows from property of
homomorphic image of nilpotent subalgebra.

Let $L_{\overline 0}$ be the Fitting's null component in the
Fitting's decomposition of the algebra $L/I$ with respect to $R(
\overline \Im).$

Suppose that $\overline \Im$ is not a Cartan subalgebra, i.e.
$\overline \Im\subset l(\overline \Im).$ By Proposition 3.1, we
have $\overline \Im \subset L_{\overline 0},$ where $L_{\overline
0}= \{ \overline x\in L/I|$ for any $\overline h\in\overline \Im$
there exists $k\in N$ such that $[\ldots[[\overline
x,\underbrace{\overline h],\overline h],\ldots}_{\mbox{$k$
times}},\overline h]=\overline 0\}.$ Hence, there exists an
element $\overline x= x+I$ such that $\overline x \in L_{\overline
0}$ and $\overline x\not\in\overline \Im$ (it means that $x\not\in
I,$ $x\not\in \Im$), i.e. $[\ldots[[\overline
x,\underbrace{\overline h],\overline h],\ldots}_{\mbox{$k$
times}},\overline h]=\overline 0$ for any $\overline h\in
\overline \Im$ and some $k\in N.$ So, there exists $x$ such that
$x\not\in\Im,$ $x\not\in I$ and for any $h\in \Im$ $[\ldots[[
x,\underbrace{ h], h],\ldots,h}_{\mbox{$k$ times}}]\in I$ for some
$k \in N.$

Thus, Lemma 3.1 is applicable to the element $x,$ i.e. $x\ne x_0,$
where $x_0\in L_{0R_b},$ and $[\ldots[[ x,\underbrace{ b],b],
\ldots, b}_{\mbox{$n$ times}}]\in I$ by some $k\in N.$ Hence
$\overline x=\overline{x_0},$ i.e. $\overline x\in \overline \Im$
which contradicts to the condition that $\overline \Im$ is not a
Cartan subalgebra.
\end{proof}

We present the example which demonstrates that the preimage by a
natural homomorphism of a regular element (Cartan subalgebra) is
not regular (Cartan subalgebra).

\emph{\bf Example 3.2.} Let the Leibniz algebra $L$ with the basis
$\{ e_1,e_2,\ldots,e_5\}$ be defined by the following
multiplication:
$$
\begin{array}{lll}
[e_2,e_1]=-e_3,& [e_1,e_2]=e_3,& [e_1,e_3]=-2e_1, \\
{[}e_3,e_1]=2e_1,& [e_3,e_2]=-2e_2,& [e_2,e_3]=2e_2, \\
{[}e_5,e_1]=e_4,& [e_4,e_2]=e_5,& [e_4,e_3]=-e_4, \\
& & [e_5,e_3]=e_5,\\ \end{array}
$$
where omitted products are equal to zero.

It is not difficult to see that $I=\{ e_4, e_5\}$ and $L/ I$ is
isomorphic to the Lie algebra $sl_2.$ Note that $\overline{e
_1}=e_1+I$ is a regular element of the algebra $L/I,$ but $e_1$ is
not regular in $L.$ Moreover, $\overline \Im =\{ \overline{
e_1}\}$ is a Cartan subalgebra of the algebra $L/I,$ but
$\Im=\{e_1,e_4,e_5\}$ is not a Cartan subalgebra of the algebra
$L.$

Similar to the case of Lie algebras case, one can easily prove
that if a nilpotent endomorphism $R_x$ corresponds to an element
$x$ of the Leibniz algebra $L$ over the field of zero
characteristic then $\exp(R_x)$ is an automorphism of the algebra
$L.$ All kinds of products of such automorphisms are said to be
\emph{invariant automorphisms}. If $\tau$ is an automorphism then
$\tau(\exp(R_x))\tau^{-1}=\exp(R_{\tau(z)}).$

\begin{thm} Let $\Im_1$ and $\Im_2$ be Cartan subalgebras of the
Leibniz algebra $L$. Then there exists an invariant automorphism
$\delta$ such that $\delta(\Im_1)=\Im_2.$ \end{thm}
\begin{proof} Let $\overline\Im$ be the image of the Cartan subalgebra $\Im$
by homomorphism $\varphi.$ Then from the theory of Lie algebras we
have the existence of a regular element $\overline a=a+I\in
\overline\Im$ such that $\overline \Im=L_{\overline 0R_{\overline
a}}.$ Take $b\in \Im$ such that $\Im=L_{0R_b}$ (existence of the
element $b$ follows from Remark 2). As $a,b\in \Im\setminus I$
then $L_{0R_b} \subseteq L_{0R_a}.$ Regularity of the element
$\overline a$ implies $L_{\overline 0 R_{\overline a}}\subseteq
L_{\overline 0 R_{\overline b}}.$ Suppose $L_{\overline 0
R_{\overline a}}\ne L_{\overline 0 R_{\overline b}}.$ Then there
exists a non-zero element $\overline x=x+I$ such that $\overline
x\in L_{\overline 0 R_{\overline b}}$ and $\overline x\not\in
R_{\overline 0 R_{\overline a}}.$ Hence, for the element $x$ we
have $[ \ldots[[x,\underbrace{b],b],\ldots, b}_{\mbox{ $k$
times}}] \in I$ for some $k\in N$ and $[ \ldots[[x,
\underbrace{a], a],\ldots, a}_{\mbox{ $s$ times}}] \not\in I$ for
any $s\in N.$

Note that $[ \ldots[[x,\underbrace{b],b],\ldots, b}_{\mbox{ $t$
times}}] \neq 0$ for any $t\in N.$ In fact, otherwise $x_0\in
L_{0R_b}$ and $x_0\in L_{0R_a}.$ But it contradicts the condition
$[ \ldots[[x,\underbrace{a],a],\ldots, a}_{\mbox{ $s$ times}}]
\not\in I$ for any $s\in N.$ So, $[ \ldots[[x,\underbrace{b],b]
,\ldots, b}_{\mbox{ $t$ times}}] \ne 0$ for any $t\in N,$ i.e. $x
\not\in \Im.$ Thus, we have obtained that for the element $x,$ the
condition $[ \ldots[[x,\underbrace{b],b],\ldots, b}_{\mbox{ $k$
times}}] \in I$ holds for some $k\in N$ and $x\ne x_0$, where
$x_0$ is the component in the decomposition of the element $x$ as
in lemma 3.1. Using this lemma, we obtain $\overline
x=\overline{x_0} \in \overline \Im = L_{\overline 0 R_{\overline
a}},$ i.e. we have the contradiction to the assumption
$L_{\overline 0 R_{\overline a}}\ne L_{\overline 0 R_{\overline
b}}.$

Thus, $\overline \Im=L_{\overline 0 R_{\overline b}}$ and
$\Im=L_{0 R_b}.$

Let $\Im _1$ and $\Im _2$ be Cartan subalgebras of the Leibniz
algebra $L.$ Then by theorem 3.2 $\overline \Im_1$ and $\overline
\Im_2$ are Cartan subalgebras of the Lie algebra $\overline L. $
Now using conjugacy of Cartan subalgebras $\overline \Im_1$ and
$\overline \Im_2$ we have that there exists an invariant
automorphism $\overline \delta$ such that $\overline \delta(
\overline \Im_1)=\overline \Im_2,$ i.e. $\delta( \overline
b_1)=\overline b_2.$ Let $\overline \delta=\exp R_{ \overline z}$,
where $\overline z=z+I$ and $z\not\in I.$ Then $(\exp R_z)(b_1)+I
=b_2+I.$

Let $c$ be an element of the Leibniz algebra $L$ such that $\exp
R_z(c)=b_2.$ Then $\exp R_z(b_1)+ I= \exp R_z(c)+I.$ It means that
$\exp R_z(b_1-c)\in I.$ Hence $b_1-c\in I$, where $b_1$ and $c$ do
not belong to the ideal $I,$ i.e. $L_{0R_c}=\Im_1,$ which
immediately implies $\exp R_z (\Im_1)=\Im_2.$
\end{proof}

Theorem 3.1 and 3.3 imply that any Cartan subalgebra of the
Leibniz algebra contains regular elements and all Cartan algebras
have the same dimension coinciding with the rank of the algebra.

{\it {\bf Acknowledgment.} The author would like to convey his
sincere thanks to professor Mukhamedov F.M. for his comprehensive
assistance and support in the completion of this paper.

Institute of Mathematics, Uzbekistan Academy of Sciences,
F.Hodjaev str.29, 700143, Tashkent (Uzbekistan), e-mail:
omirovb@mail.ru


\begin{thebibliography}{99}

\bibitem{Lod} Loday J.-L. Une version non commutative des alg$\acute
e$bres de Lie: les alg$\acute e$bres de Leibniz. L'Ens. Math., 39,
1993, p.269-293.

\bibitem{Bal} Balavoin D. Homolgy and cohomology with coefficients, of
an algebra over a quadratic operad. J. Pure Appl. Algebra, vol.
132, 1998, no. 3, p. 221-258.

\bibitem{Frab} Frabetti A. Leibniz homology of dialgebras of matrices. J. Pure Appl. Algebra,
vol. 129, 1998, p. 123-141.

\bibitem{Cas} Casas J.M. Crossed extensions of Leibniz algebras. Comm.
Algebra, vol. 27, 1999, no. 12, p. 6253-6272.

\bibitem{Cas-Pir} Casas J.M., Pirashvili T. Ten-term exact sequence of
Leibniz homology. J. Algebra, vol. 231, 2000, p. 258-264.

\bibitem{Gao} Gao Y. The second Leibniz homology group for Kac-Moody Lie
algebras, Bull. London Math. Soc., 32, 2000, p. 25-33.

\bibitem{AyupOmir1} Ayupov Sh.A., Omirov B.A. On Leibniz algebras. Algebra and
operators theory. Proceedings of the Colloquium in Tashkent 1997.
Kluwer Academic Publishers. 1998, p.1-13.


\bibitem{AlAyupOmir1} Albeverio S., Ayupov. Sh.A., Omirov B.A. On nilpotent and
simple Leibniz alge-bras. Comm. in Algebra, vol. 33, 2005, p.
1-14.

\bibitem{AyupOmir3} Ayupov Sh.A., Omirov B.A. On some classes of nilpotent
Leibniz algebras. Siberian Math. Journal, vol. 42, no. 1, 2001,
p.18-29.


\bibitem{AlAyupOmir2} Albeverio S., Ayupov. Sh.A., Omirov B.A. Cartan subalgebras
and criterion of solvability of finite dimensional Leibniz
algebras. Rheinische Friedrich-Wilhelms-Universitat Bonn,
preprint, 2004, 15pp.


\bibitem{Jac} Jacobson N. Lie algebras. Interscience Publishers, Wiley,
New York, 1962.

\end{thebibliography}
\end{document}